
 \input amssym.def
\input amssym.tex

\magnification=\magstep1  
\hsize=17,5truecm 
\vsize=25.5truecm 
\hoffset=-0.9truecm
\voffset=-0.8truecm  
\topskip=1truecm 
\footline={\tenrm\hfil\folio\hfil} 
\raggedbottom
\abovedisplayskip=3mm 
\belowdisplayskip=3mm
\abovedisplayshortskip=0mm
\belowdisplayshortskip=2mm
\normalbaselineskip=12pt  
\normalbaselines

\def\ra{\rightarrow}
\def\Z{{\Bbb Z}}
\def\Q{{\Bbb Q}}

\def\qed{\raise -2pt\hbox{\vrule\vbox to 10pt{\hrule width 4pt
                 \vfill\hrule}\vrule}}

\def\cqfd{\unskip\penalty 500\quad\qed\medbreak}

\def\simra{\buildrel\sim\over\ra}

{\bf Un th\'eor\`eme de finitude pour le groupe de Chow
des z\'ero-cycles d'un groupe alg\'ebrique
lin\'eaire  sur un corps $p$-adique}

\bigskip

Jean-Louis COLLIOT-TH\'EL\`ENE

\bigskip

{\bf Zusammenfassung} \hskip4mm  {\sl Sei $X$ eine glatte
Kompaktifizierung einer zusammenh\"angenden linearen Gruppe \"uber einem
K\"orper $k$. Die  Chowgruppe  der nulldimensionalen Zyklen von $X$ vom
Grad Null ist eine Torsionsgruppe. Wir zeigen : wenn $k$ ein $p$-adischer 
K\"orper ist, dann ist der prim-zu-$p$ Anteil dieser Gruppe
endlich. }

\bigskip

 Soit   $X$
une vari\'et\'e projective, lisse, g\'eom\'etriquement irr\'eductible
sur un corps $k$.
On note $CH_0(X)$ le groupe de Chow des z\'ero-cycles 
 modulo l'\'equivalence rationnelle, et  on note  $A_0(X) \subset CH_0(X)$
le groupe de Chow r\'eduit des z\'ero-cycles, qui est 
le sous-groupe form\'e des classes de z\'ero-cycles de degr\'e z\'ero.

Si $X$ est une vari\'et\'e rationnelle (c'est-\`a-dire birationnelle
\`a un espace projectif apr\`es extension du corps de base),
ou plus g\'en\'eralement si $X$ est une vari\'et\'e 
rationnellement connexe (au sens de Koll\'ar, Miyaoka et Mori),
alors le groupe $A_0(X)$ est un groupe de torsion, d'exposant fini. Si de
plus $k$ est un corps $p$-adique, par quoi l'on entend dans tout cet
article une extension  finie du corps ${\bf Q}_p$ des nombres
$p$-adiques, on
conjecture ([CT95], [KoSz03]) que le groupe $A_0(X)$ est fini. Le but du
pr\'esent article est d'\'etablir un cas particulier de cette conjecture.

\medskip

{\bf Th\'eor\`eme.} {\sl Soient $k$ un corps $p$-adique, $G$ un
$k$-groupe lin\'eaire connexe et $X$ une $k$-compactification lisse
de $G$, c'est-\`a-dire une $k$-vari\'et\'e projective
 lisse
contenant $G$ comme ouvert dense.
Le groupe $A_0(X)$ est la somme d'un groupe fini
et d'un groupe de torsion $p$-primaire d'exposant fini.}

\medskip

 Au paragraphe 1,
on donne une variante (Th\'eor\`eme 3) du calcul, d\^u \`a 
P. Gille [Gi97] et \`a Borovoi et Kunyavski\v{\i} [BoKu04], de la
$R$-\'equivalence sur le groupe
$G(k)$ des points
$k$-rationnels de
$G$. C'est un point essentiel pour la d\'emonstration.
Au paragraphe 2,  on donne deux lemmes sur les extensions
de corps $p$-adiques. La d\'emonstration du th\'eor\`eme
est donn\'ee au paragraphe 3. C'est la pr\'esence possible
de ramification sauvage qui emp\^eche de contr\^oler la
partie $p$-primaire du groupe $A_0(X)$.
Au paragraphe 4, on donne une minoration du groupe $A_0(X)$,
qui assure en particulier que ce groupe n'est pas toujours nul.
On donne des cas ($G$  simplement connexe ou adjoint)
o\`u $A_0(X)=0$. 
Au paragraphe 5, on rappelle ce qui est
connu sur 
la conjecture g\'en\'erale mentionn\'ee ci-dessus.

\bigskip

{\bf \S 1. R\'esolutions flasques des groupes r\'eductifs
et $R$-\'equivalence}

\medskip

Commen\c cons par quelques rappels ([CTSa77], [Vo77], [Vo98]).
Soient $k$ un corps, ${\overline k}$ une cl\^oture s\'eparable
et $\frak{g}={\rm Gal}({\overline k}/k)$.
 Un $k$-groupe de type multiplicatif  (resp. un $k$-tore) est un
$k$-groupe alg\'ebrique lin\'eaire qui, sur ${\overline k}$, se plonge
dans (resp. est isomorphe \`a) un produit de groupes multiplicatifs ${\bf
G}_m$. Pour tout entier $n>0$, on note $\mu_n \subset {\bf G}_m$
le groupe des racines $n$-i\`emes de l'unit\'e.
A tout $k$-groupe de type multiplicatif 
$M$ on associe deux modules
galoisiens de type fini (groupes ab\'eliens de type fini \'equip\'es
d'une action continue discr\`ete de $\frak{g}$),
le groupe des caract\`eres $X^*(M)={\rm Hom}(M,{\bf G}_m)$ (sur
${\overline k}$) de
$M$  et le groupe des cocaract\`eres $X_*(M)={\rm Hom}({\bf G}_m,M) $
(sur ${\overline k}$) de
$M$. Le
$k$-groupe de type multiplicatif
$M$ est un $k$-tore si et seulement si le groupe ab\'elien $X^*(M)$ 
est
sans torsion. Un $k$-tore $T$
est dit {\it quasi-trivial} si 
le module galoisien $X^*(T)$ (ou, de fa\c con \'equivalente,
$X_*(T)$) est un
$\frak{g}$-module de permutation, c'est-\`a-dire qu'il poss\`ede une base
sur ${\Z}$ 
respect\'ee par
$\frak{g}$. Un tel $k$-tore est un produit de
restrictions \`a la Weil $R_{k_i/k}{\bf G}_m$, pour certaines
extensions finies de corps $k_i/k$.
Un $k$-tore $T$ est dit {\it flasque} si pour tout sous-groupe
ouvert $\frak{h} \subset \frak{g}$, le groupe de cohomologie
$H^1(\frak{h}, X_*(T))$ est  nul. Un $k$-tore quasi-trivial est flasque.
 La cohomologie
utilis\'ee dans cet article est la cohomologie galoisienne ([Se94]).
Si $X$ est une vari\'et\'e sur un corps $k$, par quoi l'on entend
un $k$-sch\'ema s\'epar\'e de type fini, et si $K$
est un corps contenant  $k$, on note $X_K$ la $K$-vari\'et\'e
$X\times_kK$. Si $T$ est un $k$-tore quasi-trivial, resp. flasque, pour
toute extension de corps $K/k$ le $K$-tore $T_K$ est quasi-trivial, resp.
flasque.

\bigskip

Les deux \'enonc\'es rassembl\'es dans la proposition
suivante sont essentiellement dus \`a S. Endo et T. Miyata [EnMi74].
Rappelons qu'une extension galoisienne finie de corps
est dite m\'etacyclique si tout sous-groupe de Sylow
de son groupe de Galois est cyclique.

\medskip

 {\bf Proposition 1.} {\sl 
(i) Etant donn\'e un $k$-groupe de type multiplicatif $M$,
il existe une suite exacte de $k$-groupes de type multiplicatif
$$ 1 \to M \to S \to P \to 1$$
avec $S$ un $k$-tore flasque et $P$ un $k$-tore quasi-trivial.
Si $M$ est d\'eploy\'e par une extension $K/k$, on peut
choisir 
$S$ et $P$ 
d\'eploy\'es par cette extension.

(ii) Si un $k$-tore flasque $S$ est d\'eploy\'e par
une extension m\'etacyclique $K/k$, alors il existe
un $k$-tore $S_1$ tel que $S \times_k S_1$ soit
un $k$-tore quasi-trivial. En particulier $H^1(k,S)=0$.}

\medskip

{\it R\'ef\'erences.}  Pour (i), voir [CTSa87], Lemma 0.6.
Pour (ii), voir [CTSa77], Prop. 2 p.~184 ou [Vo98],
4.8, Thm. 3 p.~55. \cqfd
\bigskip

{\bf Proposition 2.} {\sl  Soient $k$ un corps de caract\'eristique
z\'ero et
$G$ un
$k$-groupe r\'eductif connexe. Il existe une suite exacte de
$k$-groupes alg\'ebriques r\'eductifs connexes
$$ 1 \to S \to H \to G \to 1$$
dans laquelle $S$ est un $k$-tore flasque,
sous-groupe central dans le groupe $H$,
et le groupe $H$ est une extension d'un $k$-tore
quasi-trivial par un $k$-groupe semi-simple simplement connexe.}

\medskip

{\it D\'emonstration.} Soit $G' \subset G$ le groupe
d\'eriv\'e de $G$. C'est un $k$-groupe semi-simple.
Soit $G^{sc}$ le rev\^etement simplement connexe de $G'$.
Soit $T \subset G$ le radical de $G$, c'est-\`a-dire la composante neutre
du centre de
$G$. C'est un $k$-tore. Soit $Q \to T$ un $k$-homomorphisme
d'un $k$-tore quasi-trivial $Q$ sur $T$ (on peut par exemple
trouver un tel homomorphisme de noyau un $k$-tore).
On dispose alors de l'homo\-mor\-phisme de $k$-groupes
alg\'ebriques 
$ G^{sc} \times_k Q \to G$ obtenu par composition.
On v\'erifie  que le noyau de cet homomorphisme est un
$k$-groupe de type multiplicatif $M$, central dans $G^{sc} \times_k Q$.
Soit
$$1 \to M \to S \to P \to 1$$ une suite exacte
de $k$-groupe de types multiplicatifs
donn\'ee par la Proposition 1 (i).
Soit  $H$ le conoyau de la fl\`eche diagonale
$M \to (G^{sc} \times_k Q) \times_k S$.
On dispose alors de la suite exacte de $k$-groupes
$$ 1 \to S \to H \to G \to 1,$$
avec $S$ central dans $H$, et d'une suite exacte
de $k$-groupes
$$ 1 \to G^{sc} \times_k Q \to  H \to P \to 1.$$
Le groupe d\'eriv\'e de $H$ est le groupe $G^{sc}$,
qui est simplement connexe, et le quotient de 
$H$ par ce $k$-sous-groupe normal est un $k$-tore 
extension du $k$-tore quasi-trivial $P$
par le 
$k$-tore
quasi-trivial $Q$. Comme toute 
extension de tels $k$-tores est scind\'ee (ce qu'on voit ais\'ement
sur les suites duales de caract\`eres), ce
quotient est isomorphe au $k$-tore quasi-trivial $P \times_k Q$.\cqfd

\bigskip

{\bf D\'efinition.} {\it Une suite exacte de $k$-groupes
alg\'ebriques r\'eductifs connexes
$$ 1 \to S \to H \to G \to 1,$$
avec $H$  extension d'un
$k$-tore quasi-trivial par un $k$-groupe semi-simple simplement connexe,
et $S$ un $k$-tore flasque central dans $H$,
est appel\'ee  une
 {\bf r\'esolution flasque} du $k$-groupe r\'eductif connexe $G$.}

\medskip

Il y a pour ces
r\'esolutions flasques des propri\'et\'es de presque unicit\'e analogues
\`a celles connues dans le cas o\`u $G$ est un  $k$-tore ([CTSa77],
[Vo77], [Vo98]). 
Ceci fera l'objet d'un expos\'e s\'epar\'e (voir l'annonce [CT04]).
Pour toute extension de corps $K/k$, la suite associ\'ee $ 1 \to S_K \to
H_K
\to G_K
\to 1$   est une r\'esolution flasque de
$G_K$.

\bigskip

Rappelons la d\'efinition de la $R$-\'equivalence
sur les $k$-points d'une vari\'et\'e alg\'ebrique $X$
d\'efinie sur un corps $k$. C'est la relation d'\'equivalence
engendr\'ee par la relation \'el\'ementaire suivante~:
deux $k$-points $A$ et $B$ de $X(k)$ sont li\'es s'il
existe un $k$-morphisme $f : U \to X$ d'un ouvert $U$ de la droite
projective ${\bf P}^1_k$ vers $X$ tel que $A$ et $B$ 
appartiennent \`a $f(U(k))$. 
Si $X$ est projective,
on peut prendre $U={\bf P}^1_k$ dans cette d\'efinition.
Si $X$ est projective, et $A$ et $B$ sont deux $k$-points 
$R$-\'equivalents, le z\'ero-cycle $A-B$ est rationnellement
\'equivalent \`a z\'ero sur $X$.

Soit $G$ un $k$-groupe alg\'ebrique. La structure de
groupe sur $G(k)$ induit une structure de groupe sur
$G(k)/R$.
On renvoie \`a [CTSa77] et [Gi97] pour plus de d\'etails.

\bigskip
On peut, de diverses fa\c cons, \'etablir la finitude
du quotient $G(k)/R$ pour $G$ un $k$-groupe lin\'eaire connexe
sur un corps $p$-adique. Mais pour \'etablir le th\'eor\`eme principal
du pr\'esent article nous aurons  besoin de la valeur pr\'ecise
du groupe $G(k)/R$, donn\'ee par l'\'enonc\'e suivant, qui est une
variante d'un r\'esultat de P. Gille ([Gi97], III.2.7) et de
Borovoi et Kunyavski\v{\i} ([BoKu04],  Thm. 4.8).

\medskip

{\bf Th\'eor\`eme 3.} {\sl Soit $k$ un corps 
$p$-adique, et soit 
$$1 \to S \to H \to G \to 1$$
une r\'esolution flasque d'un $k$-groupe r\'eductif connexe $G$.
L'application bord
$G(k) \to H^1(k,S)$ d\'eduite de cette suite
induit un isomorphisme de groupes ab\'eliens finis
$G(k)/R \simra H^1(k,S).$}

\medskip

{\it D\'emonstration.} 
Sur un corps $k$ quelconque, une suite exacte de $k$-groupes
alg\'ebriques 
$$1 \to S \to H \to G \to 1,$$ 
comme ci-dessus  induit ([Se94], I, \S 5) une
suite exacte de groupes
$$H(k) \to G(k) \to H^1(k,S), $$
et l'image de l'application $G(k) \to H^1(k,S)$
co\"{\i}ncide avec l'ensemble des \'el\'ements
de $H^1(k,S)$ d'image la classe triviale dans l'ensemble $H^1(k,H)$. 
Deux points $\alpha,\beta \in G(k)$ ont m\^eme image
dans $H^1(k,S)$ si et seulement si $\alpha.\beta^{-1} \in G(k)$
est l'image d'un \'el\'ement de $H(k)$.
Lorsque le $k$-tore $S$ est flasque, l'homo\-mor\-phisme
$\delta : G(k) \to H^1(k,S)$ passe au quotient par la $R$-\'equivalence
([CTSa77], Prop.~12 p.~198), 
on obtient un homomorphisme $G(k)/R \to H^1(k,S)$.
 
Pour \'etablir que, sur un corps $p$-adique, l'homomorphisme
$G(k)/R \to H^1(k,S)$ est un isomorphisme, il suffit de
d\'emontrer que sur un tel corps le groupe $H$ satisfait les deux
propri\'et\'es suivantes~:

(i) On a $H^1(k,H)=1$, en d'autres termes tout espace
principal homog\`ene sous $H$ (sur $k$) est trivial.

(ii) Le quotient $H(k)/R$ est r\'eduit \`a un \'el\'ement.

Le groupe $H$ s'ins\`ere dans une suite exacte
$$1 \to H' \to H \to P \to 1,$$
o\`u $H'$ (groupe d\'eriv\'e de $H$) est un $k$-groupe semi-simple
simplement connexe et o\`u $P$ est un $k$-tore quasi-trivial.
Cette suite induit une suite exacte d'ensembles point\'es
([Se94],  I, \S 5.5)
$$ H^1(k,H') \to H^1(k,H) \to H^1(k,P).$$
Le groupe $H^1(k,P)$ est trivial (th\'eor\`eme 90
de Hilbert).  Le corps $k$ \'etant $p$-adique 
et le groupe $H'$ semi-simple
simplement connexe,
l'ensemble $H^1(k,H')$ est r\'eduit 
\`a un 
\'el\'ement (th\'eor\`eme de Kneser). Ainsi $H^1(k,H)=1$, ce
qui
\'etablit le point (i).

Le groupe $H'$ \'etant simplement connexe
et le corps $k$ $p$-adique, on sait que le quotient
$H'(k)/R$ est r\'eduit \`a un \'el\'ement (Voskresenski\v{\i}, 1979,
voir [Vo98], 18.5, Thm. 1 ;  pour un \'enonc\'e  sur des corps plus
g\'en\'eraux, rassemblant les r\'esultats de nombreux autres auteurs,
voir [CTGiPa04], Thm. 4.5).  Comme le $k$-tore quasi-trivial $P$ est un
ouvert d'un espace affine, on a
$P(k)/R=1$. D'apr\`es le Th\'eor\`eme 1 de l'appendice de P. Gille  \`a
[BoKu04], dans la situation ci-dessus ($k$
corps $p$-adique, $H'$ simplement connexe), la suite exacte de $k$-groupes
alg\'ebriques $1 \to H' \to H \to P \to 1$ induit  une suite exacte de
groupes $$ H'(k)/R \to H(k)/R
\to P(k)/R.$$ On a donc $H(k)/R=1$, ce qui \'etablit le point (ii).\cqfd

\medskip

{\it Remarque.} 
La d\'emonstration du Th\'eor\`eme III.2.7 dans [Gi97]
repose sur le Lemme III.2.8, sp\'ecifique aux corps $p$-adiques, dont
la d\'emonstration est incorrecte  mais peut \^etre corrig\'ee.
Le th\'eo\-r\`eme III.2.7 de [Gi97] vaut en fait sur des corps plus
g\'en\'eraux : voir [CTGiPa04], Thm. 4.9, qui utilise le   Th\'eor\`eme 6,
p. 308 de [Gi01], ind\'ependant  de [Gi97]. La d\'emonstration du
r\'esultat de l'appendice de [BoKu04], invoqu\'ee \`a
la fin de la d\'emonstration ci-dessus, repose aussi uniquement sur
[Gi01]. Ainsi le Th\'eor\`eme 3 vaut sur tout corps $k$ satisfaisant les 
hypoth\`eses  des  Th\'eor\`emes 1.2 et 4.5 de [CTGiPa04] : le corps $k$
est de caract\'eristique nulle, de dimension cohomologique au plus 2,
sur toute extension
finie de
$k$ exposant et indice des alg\`ebres simples centrales co\"{\i}ncident,
et 
la dimension cohomologique de l'extension ab\'elienne maximale
de $k$ est au plus 1 (cette derni\`ere hypoth\`ese n'\'etant utilis\'ee
que lorsqu'il y a des facteurs de type $E_8$).

\bigskip

{\bf Proposition 4.} {\sl  Soient $k$ un corps $p$-adique, 
$G$ un $k$-groupe r\'eductif connexe,
et
$ 1 \to S \to H \to G \to 1$
une r\'esolution flasque de $G$.
Soit $K/k$ l'extension finie galoisienne qui d\'eploie
le $k$-tore $S$. Soit $F/k$ une extension finie.
 Si le compos\'e $KF/F$ est cyclique,
alors $G_F(F)/R=1$.}

\medskip

{\it D\'emonstration.}  L'\'enonc\'e r\'esulte imm\'ediatement de la
Proposition 1 (ii) et du Th\'eor\`eme 3, qui donne un isomorphisme
fonctoriel en le corps de base.\cqfd

\bigskip

{\bf Lemme 5.} {\it Soient $k$ un corps $p$-adique, $T$ un $k$-tore
et $L/k$ une extension finie de corps, de degr\'e
premier au degr\'e du corps $K$ de d\'eploiement
de $T$. Alors l'application de restriction
$H^1(k,T) \to H^1(L,T)$ est un isomorphisme.}

\medskip

{\it D\'emonstration.}
  Soit $g$ 
le groupe de Galois de l'extension $K/k$.
On a la suite exacte de restriction-inflation
$$ 0 \to H^1(g,T(K)) \to H^1(k,T) \to H^1(K,T).$$
Comme $T_K$ est un $K$-tore d\'eploy\'e, on
a $H^1(K,T)=0$ (th\'eor\`eme 90 de Hilbert).
Ainsi $H^1(k,T) = H^1(g,T(K))=H^1(g,X_*(T)\otimes_\Z K^*)$, et le groupe
$H^1(k,T)$ est annul\'e par l'ordre de $g$,
c'est-\`a-dire par le degr\'e $[K:k]$ de $K$ sur $k$.
Comme les degr\'es des extensions $K/k$ et $L/k$ sont
premiers entre eux, le compos\'e $M=KL$ de l'extension galoisienne
$K/k$ et de $L/k$ est une extension galoisienne de $L$
de degr\'e $[K:k]$, de groupe $g$. 
L'inclusion $K^* \hookrightarrow M^*$ et la norme
$M^* \to K^*$ sont $g$-\'equivariantes. Leur compos\'e
est l'\'el\'evation \`a la puissance
$[M:K]=[L:k]$.
La restriction $H^1(k,T) \to H^1(L,T)$
s'identifie \`a la fl\`eche
$H^1(g,X_*(T)\otimes_\Z K^*) \to H^1(g,X_*(T)\otimes_\Z M^*)$
induite par l'inclusion $K^* \hookrightarrow M^*$.
La composition avec l'application
$H^1(g,X_*(T)\otimes_\Z M^*) \to H^1(g,X_*(T)\otimes_\Z K^*)$
induite par la norme est la multiplication par $[L:k]$.
Ainsi le noyau de la restriction  $H^1(k,T) \to H^1(L,T)$
est-il annul\'e par $[L:k]$.
Comme  $[K:k]$ et $[L:k]$ sont premiers
entre eux, on conclut que la restriction $H^1(k,T) \to H^1(L,T)$
est injective. (Cette premi\`ere partie de la d\'emonstration vaut sur un
corps quelconque ; il s'agit, dans le pr\'esent contexte, d'une
d\'emonstration d\'etaill\'ee, requise par le rapporteur, de la formule
g\'en\'erale
${\rm Cores}_{L/k}
\circ {\rm Res}_{k/L} = [L:k]$.)

Rappelons que $X^*(T)$ d\'esigne le groupe des 
caract\`eres du $k$-tore $T$. Si $k$ est un corps
$p$-adique,
 le cup-produit
induit
une  dualit\'e parfaite de groupes finis
$H^1(k,T) \times H^1(k,X^*(T)) \to 
{\rm  Br}(k)
={\Q}/{\Z}$ ([Se94], II.
\S 5.8, Th\'eor\`eme 6). On a de m\^eme une dualit\'e parfaite de groupes
finis
$H^1(L,T) \times H^1(L,X^*(T)) \to 
{\rm Br}(L)
={\Q}/{\Z}.$
Les suites de restriction-inflation montrent que les
fl\`eches $H^1({\rm Gal}(K/k),X^*(T)) \to H^1(k,X^*(T))$ et 
$  H^1({\rm Gal}(M/L),X^*(T)) \to H^1(L,X^*(T))$ sont des isomorphismes.
L'application de restriction $H^1(k,X^*(T)) \to H^1(L,X^*(T))$
s'identifie donc \`a l'identit\'e $H^1(g,X^*(T)) = H^1(g,X^*(T))$,
c'est un isomorphisme.  En particulier ces deux groupes ont m\^eme ordre.
Ainsi la restriction $H^1(k,T) \to H^1(L,T)$ est une injection
de groupes finis de m\^eme ordre, c'est donc un isomorphisme.\cqfd

\bigskip

{\bf Proposition 6.} { \it Soient $k$ un corps $p$-adique, 
$G$ un $k$-groupe r\'eductif connexe,
et
$ 1 \to S \to H \to G \to 1$
une r\'esolution flasque de $G$.
Soit $K/k$ l'extension finie galoisienne qui d\'eploie
le $k$-tore $S$. Supposons que
le degr\'e
 $[K:k]$ soit une puissance d'un nombre premier $l$. 
Pour $F/E/k$ des extensions finies de corps avec $[F:E]$ premier
\`a $l$, l'homomorphisme de restriction
$G_E(E)/R \to G_F(F)/R$ est un isomorphisme de groupes finis.}

\medskip

{\it D\'emonstration.} 
L'\'enonc\'e r\'esulte imm\'ediatement du Lemme 5 et 
 du Th\'eor\`eme 3, qui donne un isomorphisme
fonctoriel
en le corps de base.
\cqfd

\bigskip

{\bf \S 2. Deux lemmes sur les extensions de corps $p$-adiques}

\medskip

Soit  $K/k$ une extension galoisienne
de corps $p$-adiques, de degr\'e $l^n$,
avec $(l,p)=1$.
 Soit
$E \subset K$ la sous-extension maximale non
ramifi\'ee de $K/k$. Soit $[E:k]=l^a$
et $[K:E]=l^b$ (ainsi $l^b$ est l'indice de
ramification de $K$ sur $k$).
L'extension $E/k$ est cyclique. Comme
$l$ est diff\'erent de $p$, l'extension
d'inertie $K/E$ est aussi cyclique.
L'hypoth\`ese $l\neq p$ est  ainsi indispensable dans
le lemme suivant, et c'est la raison principale pour
laquelle on n'obtient pas la finitude de la partie
$p$-primaire du groupe de Chow dans le th\'eor\`eme
principal.
\bigskip

{\bf Lemme 7.}  {\it Soient $k$ un corps $p$-adique,
$l$ un nombre premier diff\'erent de $p$,
puis $K/k$ une extension finie galoisienne
de degr\'e $l^n$.
 Soit $F/k$ une extension finie. Si
$l^n$ divise $[F:k]$, alors le compos\'e
$KF$ de $K$ et $F$ est une extension cyclique
de $F$.}

\medskip
{\it D\'emonstration.} Soit $M/k$ la sous-extension
maximale non ramifi\'ee de $F/k$. 
L'hypoth\`ese implique que l'on a au moins l'une
des deux propri\'et\'es :
$l^a$ divise $[M:k]$ ou $l^b$
divise $[F:M]$.

Si $l^a$ divise $[M:k]$,
on a $E \subset M$, et le compos\'e
de $K$ et $F$ sur $k$ est le compos\'e
de l'extension cyclique $K/E$ et
de l'extension $F/E$, c'est donc une
extension cyclique de $F$.

 Supposons 
que $l^b$ divise $[F:M]$. Montrons
que l'extension $KF/F$ est alors non ramifi\'ee 
 (ce \`a quoi l'on s'attend selon le principe : 
{\it la ramification avale la ramification}).
 Soient $\pi_k, \pi_K, \pi_F$ des
 uniformisantes
de $k, K, F$. On peut \'ecrire $\pi_k=u.\pi_K^{l^b}$
avec $u$ unit\'e dans $K$, et 
$\pi_k=v. \pi_F^{m.l^b}$, avec $v$ unit\'e
dans $F$ et $m$ entier.
Comme les extensions $K/E$ et $F/M$ sont
totalement ramifi\'ees, et que $l$ est premier
\`a $p$, les fl\`eches naturelles d'inclusion
d'unit\'es $O_E^* \subset O_K^*$ et $O_M^* \subset O_F^*$
induisent des isomorphismes 
 $O_E^*/(O_E^*)^{l^t} \simra O_K^*/(O_K^*)^{l^t}$
et $O_M^*/(O_M^*)^{l^t} \simra O_F^*/(O_F^*)^{l^t}$
pour tout entier $t \geq 1$.
  
 Quitte \`a changer d'uniformisante pour
$K$  on peut donc supposer $u \in O_E^*$.
Par ailleurs on peut \'ecrire $v=w.z^{l^b}$
avec $w \in O_M^*$ et $z \in O_F^*$.
On a donc $\pi_k=w.(z.\pi_F^m)^{l^b}$.
De ces \'equations on tire 
$$(\pi_K.(z.\pi_F^m)^{-1})^{l^b}= w.u^{-1}.$$
Comme $l$ est premier \`a $p$, et que
$w$ et $u$ sont dans des extensions non ramifi\'ees
de $k$,
 il existe
une extension non ramifi\'ee de $k$, contenant $E$,
 qui 
contient $  \pi_K.(z.\pi_F^m)^{-1}$.
 Ainsi $\pi_K$ appartient au compos\'e
de cette extension non ramifi\'ee et
de $F$, l'\'el\'ement $\pi_K$
de $KF$ appartient \`a une
extension non ramifi\'ee de $EF$.

Comme $\pi_K$ engendre $K$ sur $E$
(solution d'une \'equation d'Eisenstein),
$KF$ est une extension non ramifi\'ee de
$EF$, donc de $F$. L'extension de corps
locaux  $KF/F$
\'etant non ramifi\'ee, elle est cyclique.\cqfd

\bigskip

{\bf Lemme 8.} {\it Soient $k$ un corps $p$-adique
et $l \neq p$ un nombre premier. Supposons
que $k$ contient  les racines $l$-i\`emes de
l'unit\'e. Soit $F/k$ une extension finie de corps.  
Soit $l^n$ la plus grande puissance de $l$
divisant le degr\'e $[F:k]$. Il existe alors
une sous-extension $L/k$ de $F/k$ telle que
$[L:k]=l^n$, et donc $([F:L],l)=1$.}

\medskip

{\it D\'emonstration.} Soit $E/k$ la sous-extension
maximale non ramifi\'ee de $F/k$ et  soit $N \subset E$
 la sous-extension maximale non ramifi\'ee de $F/k$
de degr\'e  une puissance de $l$, soit $l^a$.
L'extension $E/N$ est de degr\'e premier \`a $l$.
L'inclusion  $\kappa_N \subset \kappa_E$ de leurs corps
 r\'esiduels induit, pour tout entier $t \geq 1$,
une inclusion $\kappa_N^*/\kappa_N^{*l^t}
\hookrightarrow  \kappa_E^*/\kappa_E^{*l^t}$.
L'ordre du premier quotient est le m\^eme que
celui de $\mu_{l^t}(\kappa_N)$, et celui du
deuxi\`eme quotient est le m\^eme que celui
de $\mu_{l^t}(\kappa_E)$.
L'inclusion $\mu_{l^t}(\kappa_N) \subset \mu_{l^t}(\kappa_E)$
est un isomorphisme  ; en effet, s'il existait $\zeta$ dans
le second groupe non dans le premier, la sous-extension
$\kappa_N(\zeta)/\kappa_N$ de $\kappa_E/\kappa_N$
serait de degr\'e une puissance de $l$ (car
l'hypoth\`ese sur $k$ implique que le corps
r\'esiduel de $k$ contient les racines $l$-i\`emes de 
l'unit\'e), et ceci n'est pas possible puisque
le degr\'e de $\kappa_E/\kappa_N$ est \'egal 
\`a celui de $E/N$, donc est premier \`a $l$.
Ainsi l'inclusion     $\kappa_N^*/\kappa_N^{*l^t}
\hookrightarrow  \kappa_E^*/\kappa_E^{*l^t}$
est un isomorphisme. Les corps r\'esiduels
de $E$ et $F$ co\"{\i}ncident.
En utilisant  le lemme de Hensel, on voit alors que
pour tout $t$ entier, $t \geq 1$, l'inclusion
naturelle des groupes d'unit\'es $O_N^* \subset O_F^*$
induit un isomorphisme
 $O_N^*/(O_N^*)^{l^t} \simra O_F^*/( O_F^*)^{l^t}$.
Soit $e=l^b.m$ l'indice de ramification de $F$ sur $k$, avec
$(m,l)=1$. Notant $\pi_k$ et $\pi_F$
des uniformisantes de $k$, resp. $F$, on a
  $\pi_k=u.(\pi_F)^{l^b.m}$, avec $u \in O_F^*$.
D'apr\`es ce qui pr\'ec\`ede, on peut donc \'ecrire
$v.\pi_k=w^{l^b}$, avec $v \in O_N^*$ et $w \in F^*$.
L'\'equation $X^{l^b}- v.\pi_k$ est une \'equation
d'Eisenstein sur le corps $N$, elle d\'efinit
une extension $L$ totalement ramifi\'ee de degr\'e
$l^b$ de $N$. D'apr\`es ce qui pr\'ec\`ede, 
cette \'equation a la racine $w$ dans $F$.
Ainsi $L \subset F$. On a  $[L:k]=[L:N][N:k]=l^b.l^a=l^n$.\cqfd

\medskip

{\it Remarques.} Comme me l'a fait observer P. Deligne, on peut
directement d\'eduire le Lemme~8
du fait que le groupe de Galois absolu de $k$ est
extension d'un pro-$l$-groupe par un groupe d'ordre
(profini) premier \`a $l$. Par ailleurs,
M. Brion note que lorsque l'extension finie $F/k$
est galoisienne, l'\'enonc\'e r\'esulte
(m\^eme pour $l=p$, et sans rien supposer sur les racines
de l'unit\'e) des 
th\'eor\`emes de Sylow g\'en\'eralis\'es  pour les groupes
r\'esolubles, dus \`a 
 Philip Hall (voir 
[Ha59],
Thm. 9.3.1).

\bigskip

{\bf \S 3. 
D\'emonstration du th\'eor\`eme}

\medskip

Rappelons pour m\'emoire ([Fu84], 1.4) que pour un morphisme propre
 $\pi  : Y \to Z$ entre $k$-vari\'et\'es, l'application image
directe $\pi_* : CH_0(Y) \to CH_0(Z)$ est induite par
l'application lin\'eaire sur les z\'ero-cycles $Z_0(Y) \to Z_0(Z)$ qui
envoie un point ferm\'e $M$ sur $[k(M):k(N)] N$, o\`u $N$
est le point ferm\'e de $Z$ image de $M$, et o\`u
$[k(M):k(N)]$ est le degr\'e relatif des extensions r\'esiduelles.

Soient $k$ un corps et $X$ une $k$-vari\'et\'e
projective et lisse. Soit $K/k$ une extension finie de corps.
Soit $p : X_K \to X$ la projection naturelle.
Cette projection est finie et plate.
 On dispose donc des homomorphismes
d'image r\'eciproque $p^*
: CH_0(X)
\to CH_0(X_K)$ ([Fu75], 1.2 ; [Fu84], 1.7) et d'image directe, encore
appel\'ee trace, 
$p_* : CH_0(X_K) \to CH_0(X)$ ([Fu75], 1.6 ; [Fu84], 1.4). 

Le
compos\'e
$p_* \circ p^* : CH_0(X) \to CH_0(X)$ est la multiplication
par le degr\'e $[K:k]$. C'est un cas particulier d'un
\'enonc\'e  g\'en\'eral sur les morphismes finis et plats ([Fu84],
Example 1.7.4).

En particulier, si $l$
est un nombre premier ne divisant pas
$[K:k]$, alors l'application de restriction $A_0(X) \to A_0(X_K)$
est injective sur le sous-groupe de torsion $l$-primaire de $A_0(X)$.

Le groupe de Chow r\'eduit $A_0(X)$ associ\'e
\`a une $k$-vari\'et\'e $X$ projective, lisse, g\'eom\'etri\-quement
connexe, est un invariant $k$-birationnel 
([CTCo79], Prop. 6.3 ; [Fu84], 16.1.11, o\`u l'hypoth\`ese faite
sur le corps de base est inutile). Il
satisfait
$A_0(X
\times_k Z)
\simeq A_0(X)$ si $Z$ est une $k$-vari\'et\'e projective lisse
$k$-birationnelle
\`a un espace projectif.

En caract\'eristique z\'ero, tout $k$-groupe lin\'eaire connexe $G$
est $k$-isomorphe au produit de son radical unipotent, qui est comme
$k$-vari\'et\'e un espace affine standard, par
un $k$-groupe r\'eductif connexe $G_1$. 
Pour $X$, resp. $X_1$,  une $k$-compactification lisse de $G$,
resp. de $G_1$, on a donc $A_0(X) \simeq A_0(X_1)$.
(Que de telles compactifications lisses existent 
r\'esulte du th\'eor\`eme d'Hironaka. Une d\'emonstration
plus \'econome n'est sans doute pas hors d'atteinte, mais
elle n'est pas disponible dans la litt\'erature.) 

Pour \'etablir le th\'eor\`eme, on supposera donc $G$
r\'eductif connexe.
Soit $$ 1 \to S \to H \to G \to 1$$
une r\'esolution flasque de $G$.
Soit $K/k$ l'extension finie galoisienne qui d\'eploie
le $k$-tore $S$.

Pour chaque premier $l$ divisant $[K:k]$, choisissons un $l$-sous-groupe
de Sylow $g_l$ de $g={\rm Gal}(K/k)$. Soit $k_l$ le
corps fixe de $g_l$. Comme le degr\'e de $k_l$ sur $k$
est premier \`a $l$, l'argument g\'en\'eral de trace rappel\'e 
ci-dessus assure que  
le sous-groupe de torsion
$l$-primaire de $A_0(X)$ s'injecte (par la fl\`eche de restriction
de $k$ \`a $k_l$) dans le
sous-groupe de torsion $l$-primaire de $A_0(X_{k_l})$.

Pour \'etablir le r\'esultat de finitude annonc\'e,
qui porte sur tous les corps $p$-adiques,
on peut donc supposer que le degr\'e de
 $K$ sur $k$ est une puissance d'un nombre premier $l\neq p$.
On dispose  de la r\'esolution flasque
$$ 1 \to S_K \to H_K \to G_K \to 1,$$
o\`u $S_K$ est un $K$-tore d\'eploy\'e.
Un cas  simple de la Proposition 10 ci-dessous
(qui n'utilise que le Th\'eor\`eme 3 ci-dessus) montre
alors qu'on a $A_0(X_K)=0$. L'argument de trace rappel\'e ci-dessus
montre que le groupe $A_0(X)$ est annul\'e par $[K:k]$,
qui est une puissance de $l$. (On pourrait se dispenser
de l'utilisation de la Proposition 10, en partant d'une
extension finie galoisienne $K/k$ d\'eployant le groupe $G$, et en
consid\'erant les sous-groupes de Sylow de ${\rm Gal}(K/k)$.)

Soit $\mu_l$ le groupe des racines $l$-i\`emes
de l'unit\'e. L'extension $k(\mu_l)/k$ est de degr\'e
premier \`a $l$. L'argument de trace montre 
que l'application de restriction $A_0(X) \to A_0(X_{k(\mu_l)})$
est injective. Pour \'etablir la finitude de $A_0(X)$,
il suffit d'\'etablir celle de $A_0(X_{k(\mu_l)})$. 

En r\'esum\'e, il suffit d'\'etablir la finitude de 
$A_0(X)$ lorsque  l'extension $K/k$ qui d\'eploie $S$ est de degr\'e une
puissance $l^n$ d'un nombre premier $l \neq p$
et que de plus $k$ contient $\mu_l$, deux
hypoth\`eses que nous faisons jusqu'\`a la fin de la d\'emonstration.

Pour \'etablir le r\'esultat, comme $A_0(X)$
est un groupe de torsion, il suffit de montrer
que le groupe $CH_0(X)$ est un groupe de type fini. 
Soit $M$ un point ferm\'e de $G$. Soit $F=k(M)$ le corps r\'esiduel
de $M$. Soit 
$l^t$ 
la puissance maximale de $l$ divisant
$[F:k]$.

Via l'application \'evidente $F \otimes_k F \to F$ 
envoyant $a \otimes b$ sur $ab$, 
le point ferm\'e $M$ de $G \subset X$ de
corps r\'esiduel $F$  d\'efinit  un point $F$-rationnel de $G_F \subset
X_F$, que nous noterons
$m$. L'application $p_* : CH_0(X_F) \to CH_0(X)$ envoie la classe
de $m \in X(F)$ dans $CH_0(X_F)$ sur la classe de $M$ dans $CH_0(X)$.

 Si
$n \leq t$, 
alors par
le Lemme 7 l'extension $KF/F$ est cyclique, donc d'apr\`es la
Proposition 4 le point rationnel $m$ est $R$-\'equivalent \`a
$\epsilon_F$  sur $G_F$ (on note $\epsilon \in G(k)$ l'\'el\'ement neutre
de $G$), donc aussi sur $X_F$, et les z\'ero-cycles $m$ et $\epsilon_F$
sont rationnellement
 \'equivalents sur $X_F$. En appliquant l'homomorphisme 
 $CH_0(X_F) \to CH_0(X)$ induit par la projection $X_F \to
X$,  on voit que $M$ est rationnellement
\'equivalent,  sur $X$, 
\`a un multiple de $\epsilon \in G(k)$. 

Supposons $t <n$. 
D'apr\`es le Lemme~8, il existe
une sous-extension $L/k$ de $F/k$, de degr\'e 
$l^t$,
avec $([F:L],l)=1$. Le $L$-tore $S_L$ est d\'eploy\'e par
une extension de $L$ de degr\'e une puissance de $l$.
D'apr\`es la Proposition~6, la restriction
de $L$ \`a $F$ induit un isomorphisme $G_L(L)/R \simra G_F(F)/R.$
Le groupe $G_L(L)/R$ est fini. Soient
$P_i \in G_L(L), i \in I_L,$ des repr\'esentants, en
nombre fini, 
de $G_L(L)/R$.
On dispose des projections $p_1 :X_F \to X_L$ et $p_2 : X_L \to X$,
dont le compos\'e est $p : X_F \to X$.
Ces projections induisent les applications d'image directe 
$p_{1*} : CH_0(X_F) \to CH_0(X_L)$ et $p_{2*}: CH_0(X_L) \to CH_0(X)$.
On a donc
$p_{2*}(p_{1*}(m))=p_*(m)=M \in CH_0(X)$.
Soit $P_i \in G_L(L)$ dont l'image dans $G_F(F)$
est dans la $R$-classe de $m$. Dans $CH_0(X_F)$,
on a $m=p_1^*(P_i)$.
Ainsi, dans $CH_0(X_L)$,
on a  $p_{1*}(m)=p_{1*} ( p_1^*(P_i)) = [F:L]P_i$, et dans
$CH_0(X)$, on a $$M=p_*(m)= p_{2*} ( p_{1*}(m)) =
 p_{2*}([F:L]P_i)= [F:L] 
 p_{2*}(P_i).$$

Un lemme de d\'eplacement bien connu (voir le compl\'ement ci-dessous)
assure que
le groupe $CH_0(X)$ est engendr\'e
par les points ferm\'es de l'ouvert $G \subset X$.
Le groupe $CH_0(X)$ est donc engendr\'e par $\epsilon$
et la famille des z\'ero-cycles $p_{L/k*}(P_i)$,
pour $L/k$ parcourant les extensions de $k$ de degr\'e
$l^t$ avec $t<n$
et, pour $L$ donn\'e, $i$ dans l'ensemble fini $I_L$. 
Comme il n'y a qu'un nombre fini d'extensions
de degr\'e donn\'e d'un corps $p$-adique
([Se94], chap. III, \S 4.2, p. 150/151), 
on conclut
que le groupe  $CH_0(X)$ est de type fini, et 
donc que $A_0(X)$ est fini.\cqfd

\bigskip

{\it Compl\'ement.} Dans la d\'emonstration,
nous avons utilis\'e le lemme de d\'eplacement suivant. {\it Soit $k$
un corps parfait infini.
Pour tout ouvert de Zariski non vide
$U$ d'une $k$-vari\'et\'e lisse int\`egre $V$,
tout z\'ero-cycle sur $V$ est rationnellement 
\'equivalent, sur $V$, \`a un z\'ero-cycle \`a support
dans $U$.} 
Sur requ\^ete du rapporteur, rappelons la d\'emonstration (dans le cas
quasi-projectif, on pourrait aussi invoquer [AK79], mais  le fait \`a
d\'emontrer ici est plus
\'el\'ementaire). Il suffit d'\'etablir le r\'esultat pour un point
ferm\'e $M \in V$.  Soit $F \subset V$ le ferm\'e compl\'ementaire de $U$.
Soit $d$ la dimension de $V$.  Dans l'id\'eal maximal de l'anneau local,
r\'egulier, ${\cal O}_{V,M}$  de $M$ sur $V$, il existe un \'el\'ement $g
\neq 0$  qui d\'efinit localement un ferm\'e contenant $F$.
On peut trouver 
une cha\^{\i}ne de param\`etres r\'eguliers
$f_1,\dots,f_{d-1}$, membres d'un syst\`eme de $d$ g\'en\'erateurs
de l'id\'eal maximal de  ${\cal O}_{V,M}$, telle
que l'image de
$g$ dans l'anneau local r\'egulier ${\cal O}_{V,M}/(f_1,\dots,f_{d-1})$
ne soit pas nulle. En \'ecrivant ces \'equations au voisinage de
$M$ et en prenant l'adh\'erence sch\'ematique dans $V$,
on trouve une courbe int\`egre $C$ ferm\'ee dans $V$,
passant par $M$, r\'eguli\`ere en $M$ et non contenue dans $F$.
Soit $D \to C$ la normalisation de $C$. C'est une $k$-courbe int\`egre
r\'eguli\`ere, donc lisse sur le corps parfait $k$. En particulier
$D$ est quasi-projective. Il existe un point ferm\'e $N \in D$
d'image $M$ tel qu'au voisinage de $N$, la projection 
$D \to C$ soit un isomorphisme. 
Soit
$F_1 \subset C$ le ferm\'e propre de $C$ image r\'eciproque de $F \subset
V$ par le morphisme propre $\pi  : D \to V$. 
Comme $\pi$ est propre, le  morphisme $\pi_* : Z_0(D) \to Z_0(V)$
induit sur les groupes de  z\'ero-cycles passe au quotient par
l'\'equivalence rationnelle ([Fu84], Thm. 1.4). On a $M=\pi_*(N)$.
Sur la courbe quasi-projective lisse int\`egre $D$,
le point ferm\'e $N$ est rationnellement \'equivalent
\`a un z\'ero-cycle $z$ dont le support est \'etranger \`a $F_1$
(car l'anneau semilocal de $D$ aux points de $F_1$
a un groupe de Picard trivial). Ainsi $M$
est, sur $V$, rationnellement \'equivalent 
au z\'ero-cycle $\pi_*(z)$, dont le support est dans $U$.

\vfill\eject

{\bf \S 4. Une minoration pour le groupe de Chow
r\'eduit}

\medskip

Soient $k$ un corps de caract\'eristique z\'ero, 
$G$ un $k$-groupe r\'eductif connexe  et $X$ une
$k$-compactification lisse de $G$.
Soit $1 \to S \to H \to G \to 1$ une r\'esolution flasque de $G$.
Cette r\'esolution induit un homomorphisme $G(k)/R \to H^1(k,S)$
(voir le d\'ebut de la d\'emonstration du Th\'eor\`eme 3).
Comme $S$ est flasque, il existe
un torseur ${\cal T} \to X$  sous le $k$-groupe $S$
qui \'etend le torseur sous $S$ donn\'e par $H \to G$ ([CTSa77],
Prop. 9 p. 194).
D'apr\`es ([CTSa77], Prop.
12 p.  198), ce torseur d\'efinit un homomorphisme $CH_0(X) \to H^1(k,S)$,
lequel compos\'e avec l'application naturelle compos\'ee $G(k)/R \to
X(k)/R \to CH_0(X)$ donne la fl\`eche naturelle $G(k)/R \to H^1(k,S).$
Comme le torseur ${\cal T}$ a une fibre triviale en l'\'el\'ement
neutre $\epsilon \in G(k) \subset X(k)$, on voit que l'homomorphisme
$G(k)/R \to H^1(k,S)$ peut s'\'ecrire comme le compos\'e
des applications $G(k)/R \to X(k)/R$, puis de l'application
$X(k)/R \to
A_0(X)$ qui envoie la classe d'un $k$-point $P$ sur la
classe du z\'ero-cycle $P - \epsilon$, enfin de l'homomorphisme
$A_0(X) \to H^1(k,S)$ d\'efini par le torseur  ${\cal T}$.

\medskip

{\bf Proposition 9} {\sl  Soient $k$ un corps $p$-adique,
$G$ un $k$-groupe r\'eductif connexe  et $X$ une
$k$-compactification lisse de $G$. 
Soit $1 \to S \to H \to G \to 1$ une r\'esolution flasque de $G$.
L'application naturelle $G(k)/R \to CH_0(X)$ est
injective, et $A_0(X) \to H^1(k,S)$ 
est un homomorphisme surjectif.}

\medskip

{\it D\'emonstration.}  D'apr\`es le Th\'eor\`eme 3,
l'homomorphisme $G(k)/R \to H^1(k,S)$ est une bijection.
Les consid\'erations pr\'ec\'edentes montrent alors que
 l'application compos\'ee $G(k)/R \to X(k)/R \to
A_0(X)$ induite par $P \mapsto (P - \epsilon)$
est une injection, et que l'homomorphisme $A_0(X) \to H^1(k,S)$
est une surjection. \cqfd

\medskip

{\it Remarque.} On peut de diverses fa\c cons voir que l'application
$G(k)/R \to X(k)/R$ est une bijection. Mais ceci laisse ouvert
(m\^eme pour $G$ un $k$-tore) les deux questions (\'equivalentes d'apr\`es
ce qui pr\'ec\`ede) :

(i)  L'application
$X(k)/R \to A_0(X)$ est-elle surjective ? (Tout z\'ero-cycle de degr\'e 1
sur $X$ est-il rationnellement \'equivalent \`a un point rationnel ?)

(ii)  
L'homomorphisme surjectif $A_0(X) \to H^1(k,S)$
est-il  une bijection ? 

En utilisant la Proposition 9 et le lien entre 
les r\'esolutions flasques de $G$ et les torseurs universels
\'etabli dans [CT04] (ceci utilise le th\'eor\`eme 3.2 de [BoKu04]), on
peut montrer que
l'accouplement naturel 
 $$ A_0(X) \times {\rm Br}(X)/{\rm Br}(k) \to {\rm Br}(k)={\Q}/{\Z}$$
entre le groupe $A_0(X)$ et le groupe de Brauer r\'eduit ${\rm
Br}(X)/{\rm Br}(k)$ de $X$ est 
non d\'eg\'en\'er\'e \`a droite.  
La question (ii) se reformule alors ainsi : cet accouplement est-il un
accouplement parfait de groupes finis ?

\bigskip

Nous pouvons maintenant montrer que le groupe $A_0(X)$ dont le
th\'eor\`eme principal de l'article assure la finitude (\`a la torsion
$p$-primaire pr\`es) n'est pas toujours nul.
Sur $k$ $p$-adique, 
il est facile de construire des $k$-tores flasques $S$
tels que $H^1(k,S) \neq 0$. L'exemple
le plus simple  correspond \`a une r\'esolution flasque
du $k$-tore des \'el\'ements de norme 1 dans une extension
biquadratique de $k$  (voir [CTSa77],
\S 6, Cor. 1 p.~207).

Etant donn\'e un tel $k$-tore $S$, on peut facilement
donner une suite exacte de $k$-tores alg\'e\-briques
$$ 1 \to S \to P \to G \to 1 $$
avec $P$ quasi-trivial. Si $X$ est une $k$-compactification
lisse de $G$, d'apr\`es ce qui pr\'ec\`ede, on a $A_0(X) \neq 0$. 

On peut aussi donner de tels exemples avec $G$ un $k$-groupe
semi-simple.  Soit $S/k$ comme ci-dessus.
Comme remarqu\'e par Ono (voir [Sa81], Lemme 1.7), le
th\'eor\`eme d'Artin sur les caract\`eres induits \`a partir de
groupes cycliques  assure l'existence d'un entier $m>0$,
de $k$-tores quasi-triviaux $P_1$ et $P_2$ et d'une isog\'enie
$$ 1 \to \mu \to S^m \times_k  P_1 \to P_2 \to 1.$$
Ici $\mu$ est un $k$-groupe fini commutatif.
Sur une extension finie galoisienne $K/k$, le groupe
$\mu_K$ est $K$-isomorphe \`a un produit de
groupes de racines de l'unit\'e, lesquels se plongent
dans des groupes sp\'eciaux lin\'eaires.
Le groupe $\mu$ se plonge dans le $k$-groupe fini
$R_{K/k}(\mu_K)$, et ce dernier se plonge (de fa\c con centrale)
dans un produit $G_1$ de descendus \`a la Weil de groupes
sp\'eciaux lin\'eaires. Soit $G$ le $k$-groupe
semi-simple quotient de $G_1$ par $\mu$.
Soit $H$ le quotient de $G_1 \times_k (S^m \times_k P_1)$
par l'action diagonale de $\mu$. On a d'une part une suite
exacte de $k$-groupes
$$1 \to G_1 \to H \to P_2 \to 1,$$
c'est-\`a-dire $H$ est extension d'un $k$-tore quasi-trivial
par un $k$-groupe semi-simple simplement connexe,
d'autre part
une suite exacte (centrale)
$$ 1 \to S^m \times_k P_1 \to H \to G \to 1,$$
qui est donc une r\'esolution flasque du $k$-groupe
semi-simple $G$.
Si $X$ est une $k$-compactification lisse de $G$,
le groupe $A_0(X)$ admet $H^1(k,S^m)= (H^1(k,S))^m  \neq 0$ comme
quotient.

\bigskip

Pour certains groupes $G$, le groupe
de Chow r\'eduit $A_0(X)$ est automatiquement nul.

\medskip

{\bf Proposition 10.} {\sl Soient $k$ un corps $p$-adique,
$G$ un $k$-groupe r\'eductif connexe et $X$
une $k$-compactification lisse de $G$. Dans chacun des cas suivants 

(i) $G$ est semi-simple simplement connexe,

(ii) $G$ est un groupe adjoint,

(iii) $G$ est un $k$-groupe absolument presque simple,

 (iv) il  existe
une r\'esolution flasque $1 \to S \to H \to G \to 1$
de  $G$ telle que le $k$-tore $S$ soit d\'eploy\'e par
une extension m\'etacyclique de $k$, 

 le groupe $A_0(X)$ est nul.} 

\medskip

{\it D\'emonstration.} Dans chacun des cas mentionn\'es ci-dessus,
on a $G(F)/R=1$ pour toute extension finie $F$ de $k$.
Pour les cas (i) \`a (iii), c'est \'etabli dans
[CTGiPa04],  Corollary  4.11. Dans le cas (iv), 
pour toute extension finie $F/k$, on
a $G(F)/R \simra H^1(F,S)$ (Th\'eor\`eme 3). Par ailleurs,
le $F$-tore $S_F$ est d\'eploy\'e par une extension m\'etacyclique
de $F$. On a donc $H^1(F,S)=1$ (Proposition~1).
Ainsi dans chaque cas tout point ferm\'e de $G$ est rationnellement 
\'equivalent sur $X$ \`a un multiple de $\epsilon \in G(k)$.
Par le lemme de d\'eplacement simple d\'etaill\'e \`a la fin du \S 3, ceci
suffit
\`a
\'etablir
$A_0(X)=0$.\cqfd

\medskip

{\it Remarque.} Soit $G$ un $k$-groupe r\'eductif connexe qui
admet une r\'esolution flasque $1 \to S \to H \to G \to 1$
telle que le $k$-tore $S$ soit d\'eploy\'e par
une extension $K/k$ mod\'er\'ement ramifi\'ee du corps $p$-adique $k$.
L'hypoth\`ese assure qu'un $p$-sous-groupe de Sylow de
${\rm Gal}(K/k)$ est cyclique. Soit $F$ le corps fixe d'un tel
sous-groupe. Comme l'extension $K/F$ est cyclique, la Proposition 10
assure  $A_0(X_F)=0$. Un argument de trace montre alors que le
sous-groupe de torsion $p$-primaire de $A_0(X)$ est nul.
Le th\'eor\`eme principal de cet article assure alors
que le groupe $A_0(X)$ est fini.

\bigskip

{\bf \S 5. Groupe de Chow des z\'ero-cycles sur les vari\'et\'es
rationnellement connexes}

\medskip

Soit $k$ un corps. Une $k$-vari\'et\'e  projective et lisse
$X$ est dite  {\it rationnellement connexe par cha\^{\i}nes},
si pour tout corps
alg\'ebriquement clos ${\Omega}$ contenant $k$, 
l'ensemble $X({\Omega})/R$
est r\'eduit \`a un \'el\'ement. 

Une $k$-vari\'et\'e projective et lisse
$X$ est dite {\it s\'epara\-ble\-ment rationnellement connexe}
s'il existe une extension $K/k$ de corps et
un $K$-morphisme $f : {\bf P}^1_K \to X_K$ 
tels que ``la"  d\'ecomposition de
l'image inverse $f^*T_X$ du fibr\'e tangent \`a
$X_K$ en somme directe de fibr\'es inversibles ${\cal O}_{{\bf
P}^1}(a_i)$ ne comporte que des $a_i$ strictement positifs.

Une vari\'et\'e s\'eparablement rationnellement connexe
est rationnellement connexe par cha\^{\i}\-nes. En
caract\'eristique z\'ero, les deux propri\'et\'es sont
\'equivalentes. En caract\'eristique z\'ero, on parlera
donc simplement de vari\'et\'es rationnellement connexes.

Ces vari\'et\'es ont, sur un corps alg\'ebriquement clos,
fait l'objet de nombreux travaux dans les quinze derni\`eres
ann\'ees (travaux de Koll\'ar-Miyaoka-Mori et Campana en particulier).
On se reportera \`a [Ko96] et [ArKo03].
 En caract\'eristique z\'ero, des exemples 
de telles vari\'et\'es (projectives, lisses) sont les vari\'et\'es
(g\'eom\'etriquement) rationnelles et plus g\'en\'eralement les
vari\'et\'es (g\'eom\'etriquement) unirationnelles. Un th\'eor\`eme
 de  Koll\'ar-Miyaoka-Mori et Campana assure que les
vari\'et\'es de Fano sont rationnellement connexes par cha\^{\i}nes.

\medskip

Soient $k$ un corps et $X$ une $k$-vari\'et\'e 
rationnellement connexe par cha\^{\i}nes.
Si $k$ est alg\'ebri\-quement clos,  
la trivialit\'e de $X(k)/R$ implique
$A_0(X)=0$.
Pour $k$ quelconque, un argument de trace sur les groupes de Chow montre
que  le
groupe $A_0(X)$ est un groupe de torsion.  L'\'enonc\'e suivant, plus
pr\'ecis, n'est pas dans la litt\'erature.

\bigskip

{\bf Proposition 11.} {\sl Soient $k$ un corps et $X$ une 
$k$-vari\'et\'e projective, lisse,  rationnellement
connexe par cha\^{\i}nes. Il existe un entier $n>0$ tel que pour tout
corps $F$ contenant $k$ le groupe $A_0(X_F)$ est annul\'e par $n$.}

\medskip

{\it D\'emonstration.} Un argument de trace permet de se ramener
au cas o\`u $X$ poss\`ede un point $k$-rationnel $P$.
Soit $\eta$ le point g\'en\'erique
de $X$ et $L=k(X)$ le corps des fonctions de $X$.
Soit $\Omega$ une cl\^oture alg\'ebrique de
$k(X)$.
Sur $X_{\Omega}$, les points $\eta$ et $P$ d\'efinissent
des points de $X({\Omega})$ qui sont $R$-\'equivalents, donc
rationnellement \'equivalents. Cette \'equivalence rationnelle
existe sur une extension finie de $L$. En prenant une trace, on trouve
qu'il existe un entier  $n>0$ tel que le z\'ero-cycle
$n(\eta -P_L) \in CH_0(X_L)$ soit rationnellement \'equivalent 
\`a z\'ero. Sur le produit $X \times_k X$, on trouve une \'equivalence
rationnelle entre $n(\Delta - (P \times_k X))$, o\`u $\Delta \subset X
\times_k X$ d\'esigne la diagonale, et un cycle 
support\'e dans un ferm\'e de la forme $X \times_k Y$, avec
$Y \subset X$ de codimension au moins~1. La correspondance
 $CH_0(X) \to CH_0(X)$ induite par $n(\Delta - (P \times_k X)) \in
CH^d(X\times_k X)$ (o\`u $d$ est la dimension de $X$) est
donc nulle (voir [Fu84], Chapitre  16). Par ailleurs, cette correspondance
co\"{\i}ncide avec  l'application qui \`a une classe de z\'ero-cycle
$z$ associe $n(z - {\rm deg}_k(z)P)$. Ceci montre que le groupe 
$A_0(X)$ est annul\'e par $n$. Comme l'\'equivalence rationnelle
sur $X \times_k X$ en induit une sur $X_F \times_FX_F$ pour toute
extension de corps $F/k$, on voit que l'on a aussi $nA_0(X_F)=0$.\cqfd

\medskip

Dans le reste de ce paragraphe, nous rappelons ce qui est connu sur la
conjecture suivante ([CT95] , [KoSz03]), que le th\'eor\`eme
principal du pr\'esent article \'etablit  dans un cas particulier.

\medskip

{\bf Conjecture} {\sl Soit  
 $X/k$ une vari\'et\'e  projective,
lisse,  rationnellement connexe sur un corps $p$-adique $k$. Le groupe
$A_0(X)$ est un groupe fini.}

\medskip

L'analogue de cette conjecture  sur  le corps ${\Bbb R}$ des r\'eels
est connue. Si $X({\Bbb R})=\emptyset$, alors $A_0(X)=0$.
Si $X({\Bbb R}) \neq \emptyset$ et $s \geq 1$
est le nombre de composantes connexes de l'espace topologique
$X({\Bbb R})$, alors $A_0(X) \simeq ({\Z}/2)^{s-1}$
(cas particulier d'un th\'eor\`eme de Ischebeck et l'auteur, 1981).

\medskip

En dimension 2, les vari\'et\'es (projectives, lisses)
s\'eparablement rationnellement connexes ne sont autres que les surfaces
(g\'eom\'etriquement) rationnelles. Des techniques de $K$-th\'eorie
(arguments de S. Bloch, th\'eor\`eme de Merkur'ev et Suslin),
gr\^ace auxquelles on peut contr\^oler la torsion dans le groupe
de Chow des cycles de codimension 2, 
ont permis  il y a vingt ans  de d\'emontrer  la conjecture dans ce cas
 : pour toute surface (g\'eom\'etriquement) rationnelle $X$ sur un corps
$p$-adique $k$, le groupe $A_0(X)$ est un groupe fini ([CT83]). 
De fait, la m\^eme m\'ethode donne aussi la finitude si $k$ est un corps
de nombres ([CT83])  -- alors qu'on est loin de pouvoir \'etablir
l'analogue du th\'eor\`eme principal du pr\'esent article sur un corps de
nombres.

\medskip

En dimension sup\'erieure \`a 2, 
le principal th\'eor\`eme de finitude pour $A_0(X)$ obtenu avant le
pr\'esent article concerne le cas des fibr\'es en quadriques au-dessus de
la droite projective, la m\'ethode reposant sur une r\'eduction
au cas des cycles de codimension 2 sur une surface, o\`u l'on utilise
le th\'eor\`eme de Merkur'ev et Suslin. On consultera
l'article de Parimala et Suresh
[PaSu95] pour les meilleurs r\'esultats obtenus
dans cette direction.

\medskip

Dans 
[KoSz03], 
Koll\'ar et Szab\'o \'etablissent que
si $k$ est un corps $p$-adique de corps r\'esiduel (fini) $ {\Bbb F}$ et 
$X$ est une $k$-vari\'et\'e projective,
lisse, g\'eom\'etriquement irr\'eductible, avec bonne r\'eduction
 s\'eparablement rationnellement connexe sur le corps
${\Bbb F}$, alors on a $A_0(X)=0$. (On trouve des compl\'ements utiles
dans [Ko04].)

\medskip

Sans hypoth\`ese de bonne r\'eduction,
 Koll\'ar ([Ko99]) montre que 
pour tout corps local localement compact de caract\'eristique z\'ero,
et toute $k$-vari\'et\'e projective, lisse,  
rationnellement connexe, l'ensemble $X(k)/R$ est fini (il montre que la
$R$-\'equivalence est ouverte).  On trouve des compl\'ements utiles dans
[Ko04]. Ce r\'esultat a pour cons\'equence le fait suivant.

\medskip

{\bf Proposition 12.} {\sl Soit $k$ un corps $p$-adique. Soit $X$ une
$k$-vari\'et\'e projective, lisse, 
 rationnellement connexe. Les propri\'et\'es
suivantes sont \'equivalentes~:

(i) Le groupe $A_0(X)$ est fini  ;

(ii) Il existe un entier $m \geq 1$ tel que tout
z\'ero-cycle sur $X$ de degr\'e au moins \'egal \`a $m$
est rationnellement \'equivalent \`a un z\'ero-cycle effectif.}

\medskip

{\it D\'emonstration.} 
Supposons $A_0(X)$ fini.  Choisissons des points ferm\'es 
$P_1, \dots, P_r$ de $X$ tels que le degr\'e (sur $k$) de tout
z\'ero-cycle sur $X$ est une combinaison lin\'eaire
des degr\'es (sur $k$) des points $P_1, \dots, P_r$. Soient $z_i, i\in
I,$ des z\'ero-cycles de degr\'e z\'ero, en nombre fini, repr\'esentant
les classes de $A_0(X)$.
Fixons une $k$-courbe $C \subset X$ projective, lisse,
g\'eom\'e\-tri\-quement int\`egre sur $X$ qui contient tous les points
$P_j$ et tous les points ferm\'es apparaissant dans le support des
$z_i$.  Une variante  du th\'eor\`eme de Bertini
([AK79], Thm. 1 et Thm. 7)
assure l'existence d'une telle courbe $C$ sur $X$.
Soit $g$ le genre de $C$. Soit $z$ un z\'ero-cycle sur
$X$ de degr\'e $N$. Il existe des entiers $n_j, j=1,\dots,r$
tels que le z\'ero-cycle $z-\sum_{j=1}^r n_jP_j$ soit de degr\'e nul.
Ce z\'ero-cycle est donc rationnellement \'equivalent, sur $X$,
\`a l'un des $z_i$. Ainsi $z$ est rationnellement \'equivalent,
sur $X$, au z\'ero-cycle $\sum_jn_jP_j + z_i$. Ce dernier cycle, de
degr\'e $N$, est support\'e sur la courbe projective et lisse $C$, de
genre $g$. Pour $N \geq g$, ce z\'ero-cycle est, par le th\'eor\`eme
de Riemann-Roch, rationnellement \'equivalent
sur $C$ \`a un z\'ero-cycle effectif. Ainsi
tout z\'ero-cycle $z$ sur $X$ de degr\'e $N \geq g$
est rationnellement \'equivalent \`a un z\'ero-cycle effectif.
On notera que l'argument donn\'e vaut sur tout corps,
et sous l'hypoth\`ese plus faible que le groupe $A_0(X)$
est un groupe de type fini.

Supposons maintenant (ii).
Soit $M_0 \in X$ un point ferm\'e, $d$ son degr\'e,
et $n=rd \geq m$. Pour tout z\'ero-cycle $z$ de degr\'e
z\'ero, le z\'ero-cycle $z+rM_0$ est rationnellement
\'equivalent \`a une somme $\sum_in_iP_i$, avec $n_i \geq 0$
avec $P_i$ point ferm\'e de degr\'e au plus $n$ et avec
 $n_i$  major\'e par $n$. 
Il n'y a qu'un nombre fini d'extensions d'un corps $p$-adique
de degr\'e donn\'e ([Se94], chap. III, \S 4.2, p. 150/151). Pour toute
telle extension $F/k$, d'apr\`es le th\'eor\`eme de Koll\'ar mentionn\'e
ci-dessus, l'ensemble
$X(F)/R$ est fini. Ainsi les points $P$ de $X$ de corps r\'esiduel
$F$ appartiennent \`a un nombre fini de classes dans $CH_0(X)$.
Ainsi le groupe $CH_0(X)$ est de type fini, et son sous-groupe de torsion
 $A_0(X)$ est fini.\cqfd

\medskip

{\it Remarque.} La d\'emonstration donn\'ee
au paragraphe 3 a consist\'e pr\'ecis\'ement \`a montrer que, sous
certaines hypoth\`eses, le groupe de Chow $CH_0(X)$  est engendr\'e par
des z\'ero-cycles effectifs de degr\'e born\'e. L'un des
premiers cas de finitude de
$A_0(X)$ ([CTCo79]) avait
\'et\'e obtenu  en \'etablissant l'\'enonc\'e (ii) pour les
surfaces  fibr\'ees en coniques au-dessus de la droite
projective.

\bigskip

{\bf Remerciements.}  L'essentiel de ce travail a \'et\'e r\'ealis\'e
\`a l'Institut Tata de recherche fondamentale  (T.I.F.R.)
de Mumbai, en d\'ecembre 2003 et janvier 2004. Je remercie
l'Institut 
pour son hospitalit\'e et
le Centre franco-indien pour la recherche avanc\'ee
(CEFIPRA, IFCPAR) pour son soutien.

\bigskip

{\bf  Bibliographie}

\medskip

[AK79] A. B. Altman and S. L. Kleiman, Bertini theorems for 
hypersurface sections containing a subscheme, Commun. Algebra {\bf 7}
(1979) 775-790.

[ArKo03] C. Araujo et J. Koll\'ar, Rational curves on varieties,
in  {\it Higher dimensional varieties and rational points}  (Budapest,
2001), K. B\"or\"oczky, Jr., J. Koll\'ar, T. Szamuely (eds.),
Bolyai Soc. Math. Stud. {\bf 12}, Springer, Berlin   (2003),
13--68.

[BoKu04]  M. Borovoi and B. Kunyavski\v{\i}, Arithmetical birational 
invariants of linear algebraic groups over two-dimensional geometric
fields,  with an appendix by P. Gille,  J. of Algebra  {\bf 276}
(2004) 292-339.

[CT83] J.-L. Colliot-Th\'el\`ene, Hilbert's theorem 90 for $K_2$, with 
application to the Chow groups of rational surfaces. Invent. math. {\bf
71} (1983) 1--20.

[CT95]  J.-L. Colliot-Th\'el\`ene, L'arithm\'etique des z\'ero-cycles 
(expos\'e aux Journ\'ees arithm\'e\-tiques de Bordeaux, Septembre 93),
Journal de th\'eorie des nombres de Bordeaux  {\bf 7} (1995) 51-73.

[CT04] J.-L. Colliot-Th\'el\`ene,  R\'esolutions flasques des groupes
r\'eductifs connexes, C. R. Acad. Sc. Paris, S\'er. I {\bf 339} (2004)
331-334.

[CTCo79] J.-L. Colliot-Th\'el\`ene et D. Coray, 
L'\'equivalence rationnelle sur les points ferm\'es des surfaces
rationnelles fibr\'ees en coniques,  Compositio math.  {\bf 39 }
(1979), no. 3, 301--332.  

[CTGiPa04] J.-L. Colliot-Th\'el\`ene, P. Gille et R. Parimala,
Arithmetic of linear algebraic groups over 2-dimensional
geometric fields,  Duke Math. J. {\bf 121}  (2004) 285--341.

[CTSa77] J.-L. Colliot-Th\'el\`ene et J.-J. Sansuc,
La $R$-\'equivalence sur les tores, Ann. Sc. \'Ec. Norm. Sup.
{\bf 10} (1977) 175--229.

[CTSa87]  J.-L. Colliot-Th\'el\`ene et J.-J. Sansuc,
Principal homogeneous spaces under flasque tori: Applications,
J. Algebra {\bf 106} (1987) 148--205.

[EnMi74] S. Endo et T. Miyata, On a classification of the
function fields of algebraic tori, Nagoya Math. J. {\bf 56} (1974)
85--104. Corrigendum, ibid. {\bf 59} (1979) 187--190.

[Fu75] W. Fulton, Rational equivalence on singular varieties,
Publ. Math. I.H.\'E.S. {\bf 45} (1975) 147-167.

[Fu84] W. Fulton, {\it Intersection Theory}, 
  Ergebnisse der Mathematik
 und ihrer Grenzgebiete, 3. Folge, Bd. 2. 
Springer, 1984 
(2. Auflage, 1998).

[Gi97] P. Gille, La R-\'equivalence sur les groupes alg\'ebriques
r\'eductifs d\'efinis sur un corps global,
Publications Math\'ematiques de l'I.H.\'E.S. {\bf 86} (1997) 
199--235.

[Gi01] P. Gille,  Cohomologie galoisienne des groupes quasi-d\'eploy\'es
sur des corps de dimension cohomologique $\leq 2$,
Compositio mathematica {\bf 125} (2001) 283--325.

[Ha59] Marshall Hall, Jr.,  {\it The theory of groups},
The Macmillan Co., New York, N.Y., 1959. Reprint, Chelsea, New York, N.
Y., 1976.

[Ko96] J. Koll\'ar, {\it Rational curves on algebraic varieties},
Ergebnisse der Mathematik
 und ihrer Grenzgebiete, 3. Folge, Bd. 32,
Springer, 1996.

[Ko99] J. Koll\'ar, Rationally connected varieties over local fields,
Annals of Math. {\bf 150} (1999) 357--367.

[Ko04] J. Koll\'ar, Specialization of zero-cycles, 
Publications of RIMS (2004), Specialization of zero cycles. 
Publ. Res. Inst. Math. Sci.  {\bf 40}  (2004),  no. 3, 689--708.

[KoSz03] J. Koll\'ar et E. Szab\'o, Rationally connected
varieties over finite fields, Duke Math. J. {\bf 120} (2003) 251--267.

[PaSu95] R. Parimala and V. Suresh, Zero-cycles on
quadric fibrations: finiteness theorems and the cycle map. Invent. math. 
{\bf 122} (1995), 83--117.  Erratum, ibid.  {\bf 123} 
(1996),   611.

[Sa81]  J.-J. Sansuc, Groupe de Brauer et arithm\'etique des
groupes alg\'ebriques lin\'eaires, J.~reine angew. Math.
(Crelle) {\bf 327} (1981) 12--80.

[Se94] J-P. Serre, {\it Cohomologie galoisienne}, Cinqui\`eme \'edition,
 Lecture Notes in Mathematics {\bf 5}, Sprin\-ger, Berlin, 1994.

[Vo77] V. E. Voskresenski\v{\i}, 
{\it Algebraicheskie tory},
Izdat. ``Nauka'', Moscou, 1977.

[Vo98] V. E. Voskresenski\v{\i}, {\it Algebraic groups and their
birational invariants}, Transl. Mathe\-ma\-tical Monographs {\bf 179},
Amer. Math. Soc., 1998.

\vskip2cm

J.-L. Colliot-Th\'el\`ene,

Centre National de la Recherche Scientifique, 

Unit\'e mixte de recherche 8628,

Math\'ematiques,

B\^atiment 425,

Universit\'e Paris-Sud,

F-91405 Orsay

France

courriel : colliot@math.u-psud.fr
 \bye